\documentclass[12pt]{article}
\usepackage{amsmath}
\newcommand{\bee}{\begin{eqnarray*}}
\newcommand{\ene}{\end{eqnarray*}}
\newcommand{\beeq}{\begin{equation}}
\newcommand{\eneq}{\end{equation}}
\newtheorem{lem}{Lemma}[section]
\newcommand{\bel}{\begin{lem}}
\newcommand{\enl}{\end{lem}}
\newtheorem{exap}{Example}[section]
\newcommand{\beex}{\begin{exap}}
\newcommand{\enex}{\end{exap}}
\newtheorem{theo}{Theorem}[section]
\newcommand{\beth}{\begin{theo}}
\newcommand{\enth}{\end{theo}}
\newtheorem{prop}{Proposition}[section]
\newcommand{\bep}{\begin{prop}}
\newcommand{\enp}{\end{prop}}
\newtheorem{cor}{Corollary}[section]
\newcommand{\bec}{\begin{cor}}
\newcommand{\enc}{\end{cor}}
\newtheorem{defi}{Definition}[section]
\newcommand{\bef}{\begin{defi}}
\newcommand{\enf}{\end{defi}}
\newtheorem{rem}{Remark}[section]
\newcommand{\ber}{\begin{rem}}
\newcommand{\enr}{\end{rem}}

\usepackage{amsmath}
\usepackage{graphicx}
\usepackage{latexsym}

\begin{document}
 
\begin{center}
{\bf 
Concerns on Monotonic Imbalance Bounding Matching Methods\footnote{Appeared as online supplement to the {\em JASA} paper by
Iacus, King and Porro (2011) with response from the authors.}
}\\
by\\
Yannis G. Yatracos\\
Cyprus U. of Technology\\
{\em e-mail:} yannis.yatracos@cut.ac.cy
\end{center}

\begin{center} \vspace{0.05in} {\large Summary} \end{center}
\vspace{0.03in}
\begin{center}
\parbox{4.8in}
{ \quad \em \  Concerns are  expressed for
 the
Monotonic Imbalance Bounding (MIB) 
property (Iacus et al. 2011) and for  MIB  matching
because
 i) the definition of the  MIB property  leads to inconsistencies and the  nature of the imbalance  measure is not clearly defined,
 ii)  MIB  property does not generalize
  Equal Percent Bias Reducing (EPBR)
 property,
 iii) MIB matching does not provide  statistical information available with  EPBR
 matching.

}
\end{center}


 Imbalance Bounding (IB) matching 
is examined but the findings  and  the comments  remain  valid for MIB matching.
Familiarity of the readers with Iacus, King and Porro (2011, hereafter IKP 2011)  and the notation therein is assumed.
 \\

{\em On the definition of the  IB  property} \\

We 
use the IB-Definition
obtained from the authors in a recent communication.\\

{\bf IB-Definition:}  Let  $f$ be any measurable function and $D(\cdot,\cdot)$ any {\em measure of imbalance}
 that can be bounded by a scalar.
Assume  (a) fixed sizes of the random samples $n_T, \  n_C,$ (b) fixed distributions of ${\bf X}, \  P_T$  for the treated population
 and $P_C$
for the control population, 
(c)  a fixed matching method is used. 
If {\em {for a given value of $\delta$}} we obtain matched samples of sizes $m_T$ and $m_C$  such 
that 
\begin{equation}
\label{eq:IB}
D[f({\cal X}_{m_T}),f({\cal X}_{m_C}]\le \delta,
\end{equation}
then we have the property IB; ${\cal X}_{m_T}, \ {\cal X}_{m_C}$ are, respectively, the matched-treated and matched-control data.\\

 Since  IB property (\ref{eq:IB})
 is dependent on a $\delta$-value determined in advance (ex-ante), 
 the following situations will occur:\\

{\em a)}  for  fixed  treatment and  control  populations  and 
 a statistician with  ex-ante
$\delta=\delta_1$  the matching   method has not  the  IB property,  but for another statistician with ex-ante $\delta=\delta_2>\delta_1$ and $\delta_2$ large enough 
the same matching  method has the IB property. 
Consequently,
the  two statisticians 
 will be in disagreement on whether the  
matching method has the IB property  or not, leading to an inconsistency. Thus,  the class of IB matching methods is not well defined.\\

 {\em b)} For readers inclined to  justify the inconsistency  in {\em a)} from
 the subjective choice of different $\delta$-values
by  the two statisticians,
consider  a fixed matching method, one statistician 
and two sets of treatment and control populations
 with distributions $P_{T_i}$ and  
$P_{C_i}$ and bounds $\delta_i,  \ i=1,2.$ This  statistician may find that the matching method satisfies 
 IB property  (\ref{eq:IB}) for the distributions and the  bound $(P_{T_1}, P_{C_1}, \delta_1)$  but  not 
for the distributions and the bound $(P_{T_2}, P_{C_2}, \delta_2).$  Does the matching method have the IB property in this situation?
This will  hold if the IB definition is population dependent
but it is not  the case
since, according to the authors,  no assumptions are needed on the populations' distributions for
 IB property to hold
(IKP 2011, p. 346, section 2.2, lines 2-4).  \\


Looking at  (\ref{eq:IB})  any graduate student in statistics would ask ``What is the probability that (\ref{eq:IB}) holds?'' given
that  $f$ is a measurable function not necessarily constant. If this probability is equal to 1,  questions
will arise concerning the applicability of the method for all populations' distributions, suggesting
that IB-definition
 is population
dependent.  If this probability is less than 1, 
IB definition is data dependent and therefore, for a fixed
$\delta$-value, 
the
matching method may  have the IB property for one data set but 
this may not hold for a different data set from the same population. \\

Note that in  IB definition  (IKP 2011, p. 347), $D(x,y)$ denotes a distance between $x$ and $y$ but in the examples 
following this definition 
$D(x,y)=E(x-y)$ and $D(x,y)=|x|$ are not distances; $E$ denotes expected value. In  IB definition (\ref{eq:IB})
``D is any measure of imbalance'' but no precise definition of what this means  is available.
 There are no guidelines for the choice of the $\delta$-value
and a natural approach for its selection
 presented below makes IB definition data dependent.\\


{\em IB  and EPBR matchings-Does IB matching generalize  EPBR matching?} \\

 Our  main argument 
against the claim in IKP (2011)  that IB property  generalizes Rubin's  EPBR property (Rubin 1976)
is
that  IB loses the EPBR property of affine invariance with respect to linear combination of means. 
In a recent communication the authors provided the arguments that follow,  in order to show 
that  IB property is ``a mathematical generalization of EPBR property.''
Their  motivation for 
the
IB definition is 
EPBR definition
\begin{equation}
\label{eq:EPBR}
\mu_{T^*}-\mu_{C^*}=\gamma (\mu_{T}-\mu_{C}), \ 0<\gamma<1,
\end{equation}
 i.e. 
 the expected value of the difference of matched samples means, 
 $\mu_{T^*}-\mu_{C^*},$   is a proportion $\gamma$ of the  
expected value of the  difference of random samples  means,  $\mu_{T}-\mu_{C}.$  
For elements $x, y$ let $D(x,y)=x-y$  and for  a random vector $A$  set $f(A)=E(A); \ E(A)$ denotes the expected value of $A.$
Then,  EPBR property (\ref{eq:EPBR}) is rewritten
in the IB-like notation
\begin{equation}
\label{eq:motivMIB}
D[f({\cal X}_{m_T}),f({\cal X}_{m_C})]=\delta,
\end{equation} 
with
\begin{equation}
\label{eq:Ddelta}
D[f({\cal X}_{m_T}),f({\cal X}_{m_C}]=\mu_{T^*}-\mu_{C^*}, \ \hspace{3ex}  \delta=\gamma (\mu_{T}-\mu_{C});
\end{equation}
${\cal X}_{m_T}, \ {\cal X}_{m_C}$
denote matched sub-samples and 
$f({\cal X}_{m_T}), \ f({\cal X}_{m_C})$ denote  the expectations of the  matched sample means.
Finally, the equality sign in (\ref{eq:motivMIB}) is replaced by ``$\le$'' and 
the authors' conclusion is that
``In this way, we have shown that IB is a direct mathematical generalization of EPBR.''
However,
$E(A)$  is a functional of the cumulative distribution  of $A$ and this also holds for the expectations' differences  
in (\ref{eq:EPBR})
{\em but}  is  not
reflected in  (\ref{eq:motivMIB}) and the IB definition  (\ref{eq:IB}) which only  involve  measurable functions of the data. Thus, the authors' arguments 
do not show  that
IB property is mathematical generalization of  EPBR property.
\\

Irrespective of the last  sentence, using the authors' motivation we examine whether statistical information other than affine invariance is lost 
with IB matching  methods.
 Going one step further from $\delta$'s  definition  in (\ref{eq:Ddelta})  we obtain from the EPBR property 
\begin{equation}
\label{eq:Ddeltastep}
 \delta=\gamma (\mu_{T}-\mu_{C})=\gamma D[f({\cal X}_{n_T}),f({\cal X}_{n_C})],
\end{equation}
${\cal X}_{n_T},{\cal X}_{n_C}$  are  random samples.
It may occur, for example, that  a  practitioner  uses an IB matching method
with (small) $\delta$-value $10^{-4},$
but the value of  $D[f({\cal X}_{n_T}),f({\cal X}_{n_C})]$   is  $10^{-6}.$
Equation (\ref{eq:Ddeltastep})  suggests  using $ D[f({\cal X}_{n_T}),f({\cal X}_{n_C})]$
to determine an appropriate
$\delta$-value 
but  this will  introduce random sample dependence in the IB definition
unless $f$ is constant.\\

Unlike the IB property, EPBR property (\ref{eq:EPBR}) provides via $\gamma$  information  on the improvement of  expected  matched  means' difference 
compared to the
expected  random means' difference.
In the IB definition a subjective $\delta$-value is used, it is not clear in IKP (2011) what this value should be 
and  there is no comparison with  the $D$-value obtained via $f$  for  random samples. 
Moreover, with the EPBR property both sides in (\ref{eq:EPBR}) have the same sign. This information  is also  lost with  a
distance $D$  in the IB  definition.
Thus,  there is loss
  of statistical information when  using IB matching methods instead of EPBR matching methods.\\

EPBR property   (\ref{eq:EPBR}) is  clearly  moments' property and mild moment  conditions are provided   in Yatracos (2013) for  EPBR property to hold for a class of matching methods,  
thus relaxing the criticism that EPBR holds only  under restricted distributional assumptions  (IKP 2011, p. 346).\\


 
The  concerns presented  for  IB  matching methods  hold  also for  MIB matching methods  (IKP  2011, p. 347)  for which (\ref{eq:IB}) holds
with  data ${\cal X}_{m_T(\pi)},\ {\cal X}_{m_C(\pi)}$ and upper bound $\gamma_{f,D} (\pi)$  (instead of $\delta$) depending on a tuning parameter $\pi;$ 
$\gamma _{f,D}(\cdot )$ is
monotonically increasing in $\pi.$\\


  

\begin{center}
{\bf References}
\end{center}

 .

Iacus, S. M.,  King, G.  and Porro, G. (2011) Multivariate Matching Methods that are Monotonic Imbalance Bounding (MIB)
{\em JASA},  {\bf 106}, p. 345-361.

Rubin, D. B. (1976) Multivariate matching methods that are Equal Percent Bias Reducing, II: Maximums on bias reduction for fixed sample sizes. {\em Biometrics}, {\bf 32}, 121-132.

Yatracos, Y. G. (2013) Equal Percent Bias Reduction and Variance  Proportionate Modifying Properties  with 
Mean-Covariance Preserving Matching. {\em Ann. Inst. Stat. Math.} {\bf 65, 1,} 69-87.


\end{document}